\documentclass[11pt]{article}

\newread\epsffilein    
\newif\ifepsffileok    
\newif\ifepsfbbfound   
\newif\ifepsfverbose   
\newif\ifepsfdraft     
\newdimen\epsfxsize    
\newdimen\epsfysize    
\newdimen\epsftsize    
\newdimen\epsfrsize    
\newdimen\epsftmp      
\newdimen\pspoints     
\def\epsfbox#1{\global\def\epsfllx{72}\global\def\epsflly{72}%
   \global\def\epsfurx{540}\global\def\epsfury{720}%
   \def\lbracket{[}\def\testit{#1}\ifx\testit\lbracket
   \let\next=\epsfgetlitbb\else\let\next=\epsfnormal\fi\next{#1}}%
\def\epsfgetlitbb#1#2 #3 #4 #5]#6{\epsfgrab #2 #3 #4 #5 .\\%
   \epsfsetgraph{#6}}%
\def\epsfnormal#1{\epsfgetbb{#1}\epsfsetgraph{#1}}%
\def\epsfgetbb#1{%
\openin\epsffilein=#1
\ifeof\epsffilein\errmessage{I couldn't open #1, will ignore it}\else
   {\epsffileoktrue \chardef\other=12
    \def\do##1{\catcode`##1=\other}\dospecials \catcode`\ =10
    \loop
       \read\epsffilein to \epsffileline
       \ifeof\epsffilein\epsffileokfalse\else
          \expandafter\epsfaux\epsffileline:. \\%
       \fi
   \ifepsffileok\repeat
   \ifepsfbbfound\else
    \ifepsfverbose\message{No bounding box comment in #1; using defaults}\fi\fi
   }\closein\epsffilein\fi}%
\def\epsfclipoff{\def\epsfclipstring{\ifepsfdraft\space clip\fi}}%
\epsfclipoff
\def\epsfsetgraph#1{%
   \epsfrsize=\epsfury\pspoints
   \advance\epsfrsize by-\epsflly\pspoints
   \epsftsize=\epsfurx\pspoints
   \advance\epsftsize by-\epsfllx\pspoints
   \epsfxsize\epsfsize\epsftsize\epsfrsize
   \ifnum\epsfxsize=0 \ifnum\epsfysize=0
      \epsfxsize=\epsftsize \epsfysize=\epsfrsize
      \epsfrsize=0pt
     \else\epsftmp=\epsftsize \divide\epsftmp\epsfrsize
       \epsfxsize=\epsfysize \multiply\epsfxsize\epsftmp
       \multiply\epsftmp\epsfrsize \advance\epsftsize-\epsftmp
       \epsftmp=\epsfysize
       \loop \advance\epsftsize\epsftsize \divide\epsftmp 2
       \ifnum\epsftmp>0
          \ifnum\epsftsize<\epsfrsize\else
             \advance\epsftsize-\epsfrsize \advance\epsfxsize\epsftmp \fi
       \repeat
       \epsfrsize=0pt
     \fi
   \else \ifnum\epsfysize=0
     \epsftmp=\epsfrsize \divide\epsftmp\epsftsize
     \epsfysize=\epsfxsize \multiply\epsfysize\epsftmp   
     \multiply\epsftmp\epsftsize \advance\epsfrsize-\epsftmp
     \epsftmp=\epsfxsize
     \loop \advance\epsfrsize\epsfrsize \divide\epsftmp 2
     \ifnum\epsftmp>0
        \ifnum\epsfrsize<\epsftsize\else
           \advance\epsfrsize-\epsftsize \advance\epsfysize\epsftmp \fi
     \repeat
     \epsfrsize=0pt
    \else
     \epsfrsize=\epsfysize
    \fi
   \fi
   \ifepsfverbose\message{#1: width=\the\epsfxsize, height=\the\epsfysize}\fi
   \epsftmp=10\epsfxsize \divide\epsftmp\pspoints
   \vbox to\epsfysize{\vfil\hbox to\epsfxsize{%
      \ifnum\epsfrsize=0\relax
        \includegraphics{\ifepsfdraft}%
      \else
        \epsfrsize=10\epsfysize \divide\epsfrsize\pspoints
        \includegraphics{\ifepsfdraft}%
      \fi
      \hfil}}%
\global\epsfxsize=0pt\global\epsfysize=0pt}%
{\catcode`\%=12 \global\let\epsfpercent=
\long\def\epsfaux#1#2:#3\\{\ifx#1\epsfpercent
   \def\testit{#2}\ifx\testit\epsfbblit
      \epsfgrab #3 . . . \\%
      \epsffileokfalse
      \global\epsfbbfoundtrue
   \fi\else\ifx#1\par\else\epsffileokfalse\fi\fi}%
\def\epsfempty{}%
\def\epsfgrab #1 #2 #3 #4 #5\\{%
\global\def\epsfllx{#1}\ifx\epsfllx\epsfempty
      \epsfgrab #2 #3 #4 #5 .\\\else
   \global\def\epsflly{#2}%
   \global\def\epsfurx{#3}\global\def\epsfury{#4}\fi}%
\def\epsfsize#1#2{\epsfxsize}
\newtheorem{theorem}{Theorem}[section]
\newtheorem{lemma}[theorem]{Lemma}

\newtheorem{proposition}[theorem]{Proposition}

\newtheorem{definition}[theorem]{Definition}

\newcommand{\qed}{\hbox{\rule{6pt}{6pt}}}

\newcommand{\cmapright}[2]{
\smash{\mathop{
\hbox to 1cm{\rightarrowfill}}\limits^{#1}_{#2}}}

\newcommand{\cmapleft}[2]{
\smash{\mathop{
\hbox to 1cm{\leftarrowfill}}\limits^{#1}_{#2}}}

\begin{document}

\title{Bordism of Unoriented Surfaces in $4$-Space}

\author{
J. Scott Carter \\
University of South Alabama \\
Mobile, AL 36688 \\ carter@mathstat.usouthal.edu
  \and
Seiichi Kamada \\
Osaka City University \\
Osaka 558-8585, JAPAN\\
kamada@sci.osaka-cu.ac.jp
\and
Masahico Saito \\
University of South Florida \\
Tampa, FL 33620 \\ saito@math.usf.edu
\and
Shin Satoh \\
Kyoto University \\
Kyoto 606-8502, JAPAN\\
satoh@kurims.kyoto-u.ac.jp
}
\maketitle


\begin{abstract}
The group of bordism classes of
unoriented surfaces
in $4$-space 
is  determined.
The bordism classes are  characterized 
by normal Euler numbers, double linking numbers,
and triple linking numbers.
\end{abstract}


\section{Introduction}

B. J. Sanderson \cite{Sanderson1, Sanderson2}
studied 
the group $L_{m,n}$
of bordism classes of
`oriented' closed $(m-2)$-manifolds
of $n$ components in ${\bf R}^m$.
He showed that
$L_{m,n}$ is isomorphic to the homotopy group
$\pi_m\left(\bigvee_{i=1}^{n-1} S^2\right)$;
in particular,
the bordism group $L_{m,n}$ for $m=4$ is given
as follows.

\begin{theorem}[Sanderson]
\begin{eqnarray*}
L_{4,n}\cong
(\underbrace{{\bf Z}_2\oplus\dots
\oplus{\bf Z}_2}_{\frac{n(n-1)}{2}})
\oplus
(\underbrace{{\bf Z} \oplus\dots
\oplus{\bf Z}}_{\frac{n(n-1)(n-2)}{3}}).
\end{eqnarray*}
In particular, we have
$L_{4,1}\cong\{0\}$,
$L_{4,2}\cong{\bf Z}_2$, and
$L_{4,3}\cong{\bf Z}_2^3\oplus{\bf Z}^2$.
\end{theorem}

Similarly,
there is a group of
bordism classes of `unoriented' closed $(m-2)$-manifolds
of $n$ components in ${\bf R}^m$.
We denote the group by $UL_{m,n}$.
The aim of this paper is
to determine 
 the bordism group $UL_{m,n}$ for $m=4$ 
via purely geometric techniques.

An {\it $n$-component surface-link}  $F=F_1\cup\dots\cup F_n$
is a disjoint union of  closed surfaces embedded in ${\bf R}^4$
(smoothly, or PL and locally flatly) such that each component 
$F_i$, $i=1, \ldots, n$, is an 
orientable or non-orientable embedded surface,
 which may or may not be connected.
Two $n$-component surface-links $F$ and $F'$ are
{\it unorientedly bordant}
if there is a compact $3$-manifold
$W=\cup_{i=1}^n W_i$ properly embedded in
${\bf R}^4\times[0,1]$
such that
$\partial W_i=F_i\times\{0\}\cup F_i'\times\{1\}$
for $i=1,\dots,n$.
In this paper,
$F\simeq_B  F'$ means that $F$ and $F'$ are
unorientedly bordant, and
$F\cong_A F'$ means that they are ambient isotopic
in ${\bf R}^4$.
The unoriented bordism classes of
$n$-component surface-links form
an abelian group $UL_{4,n}$ such that
the sum $[F]+[F']$ is defined
to be the class $[F\amalg F']$ of
the split union $F\amalg F'$.
The identity is represented by the empty $F=\emptyset$
and the inverse $-[F]$ is
represented by the mirror image of $F$.
The following is our main theorem.


\begin{theorem} \label{thm:mainA}
\begin{eqnarray*}
UL_{4,n}\cong
(\underbrace{{\bf Z}\oplus\dots
\oplus{\bf Z}}_{n})
\oplus
(\underbrace{{\bf Z}_4\oplus\dots
\oplus{\bf Z}_4}_{\frac{n(n-1)}{2}})
\oplus
(\underbrace{{\bf Z}_2 \oplus\dots
\oplus{\bf Z}_2}_{\frac{n(n-1)(n-2)}{3}}).
\end{eqnarray*}
In particular,
we have
$UL_{4,1}\cong{\bf Z}$,
$UL_{4,2}\cong{\bf Z}^2\oplus{\bf Z}_4$,
and
$UL_{4,3}\cong{\bf Z}^3
\oplus{\bf Z}_4^3
\oplus{\bf Z}_2^2$.

\end{theorem}


This paper is organized as follows.
In \S$2$,
we give definitions of
three kinds of unoriented bordism invariants ---
normal Euler numbers, double linkings,
and triple linkings ---
by the projection method.
In \S$3$,
we study $1$-component surface-links.
In \S$4$,
we introduce a family of surface-links
which are called {\it necklaces}.
\S$5$ is devoted to study
crossing changes
that produce necklaces.
In \S$6$, we prove Theorem $1.2$.


\section{Unoriented Bordism Invariants}

All of our bordism invariants will be defined using the diagram of a knotted or linked surface. We begin by recalling this notion.
Consider 
a surface-link $F$.
We may assume that
the restriction
$\pi|_F:F\rightarrow{\bf R}^3$
of a projection
$\pi:{\bf R}^4\rightarrow{\bf R}^3$
is a generic map,
that is, the singularity set of
the image $\pi(F)\subset{\bf R}^3$ consists of
double points and isolated branch/triple points.
See Figure 1.
The closure of the self-intersection set on $\pi(F)$
is regarded as
a union of immersed arcs and loops,
which we call double curves.
Branch points (or Whitney umbrella points) occur at
 the end of the double curves, 
and triple points occur when double curves intersect.

\begin{figure}[htb]
\begin{center}
\mbox{
\epsfxsize=3in
\epsfbox{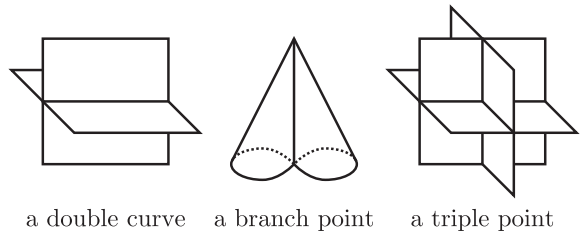} 
}
\end{center}
\caption{Generic intersections of surfaces in 3-space}
\end{figure}

A {\it surface diagram} of $F$
is the image $\pi(F)$ equipped with
over/under-information along each double curve
with respect to the projection direction.
To indicate such over/under-information,
we remove a neighborhood of
a double curve on the sheet
(the {\it lower sheet})
which lies lower than the other sheet
(the {\it upper sheet}).
See Figure 2. 
Notice that the removal of this neighborhood is merely a convention in depicting illustrations. In particular, we still speak of ``double curves'' and triple points, and we locally parametrize the surfaces using immersions.

There are seven kinds of local moves on surface diagrams,
called {\it Roseman moves},
analogues of Reidemeister moves for classical knots, 
which are sufficient to 
relate diagrams of 
ambient isotopic surface-links
(cf. \cite{CS:Reidemeister, CS:book, Roseman}).
Specifically, two diagrams represent ambiently isotopic surface-links  if and only if one can be obtained from the other by means a finite sequence of moves taken from the list of Roseman moves. We remark that the Roseman moves are used here only in relation their affect on the cobordism invariants.

\begin{figure}[htb]
\begin{center}
\mbox{
\epsfxsize=4in
\epsfbox{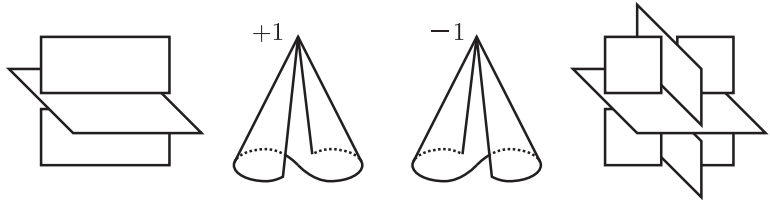} 
}
\end{center}
\caption{The broken surface diagrams at intersection points}
\end{figure}

For a surface-knot $K$
(that is a connected closed surface
embedded in ${\bf R}^4$),
H. Whitney defined the {\it normal Euler number} $e(K)$
of $K$ 
to be the 
Euler number
of a tubular neighborhood of $K$ in ${\bf R}^4$ 
considered as a $2$-plane bundle 
(cf. \cite{Massey, Whitney}).
It is known that

\begin{itemize}
\setlength{\itemsep}{-3pt}
\item[{\rm (i)}]
$e(K)=0$ if $K$ is orientable,

\item[{\rm (ii)}] (Whitney's congruence)
$e(K)\equiv 2\chi(K)$ (mod $4$), and

\item[{\rm (iii)}] (Whitney-Massey 
theorem)
$|e(K)|\leq 4-2\chi(K)$,
\end{itemize}

\noindent
where $\chi(K)$ denotes
the Euler characteristic of $K$
(cf. \cite{Kamada, Massey}).
For a $1$-component surface-link $F$,
we define the normal Euler number $e(F)$ of $F$ 
to be 
the sum of $e(K)$
for the connected components $K$ of $F$.
The normal Euler number
$e(F)$ can be calculated
by use of a projection of $F$ in ${\bf R}^3$;
it is equal to the number of positive type branch points
(see Figure~$2$) minus that of negative type ones,
\cite{CS:branch}.

Let $F=F_1\cup F_2$
be a $2$-component surface-link
and $D$ a surface diagram of $F$.
A double curve of $D$ is said to be
{\it of type $(i,j)$}
if the upper sheet 
belongs to 
$F_i$
and the lower 
belongs to 
$F_j$,
where $i,j\in\{1,2\}$.
If a double curve is an immersed arc,
then its endpoints are branch points, and hence,
the type is $(1,1)$ or $(2,2)$.
Let $C=c_1\cup\dots\cup c_m$
be the set of the double curves of type $(1,2)$
on the surface diagram $D$.
Each double curve $c_i$ is
an immersed loop in ${\bf R}^3$.

We take
a $2$-disk $B^2$ and a union of intervals $X$
in ${\bf R}^2$ as follows:

\begin{itemize}
\setlength{\itemsep}{-3pt}
\item[{\rm (i)}]
$B^2=\{(x,y)|x^2+y^2\leq 1\}$, and

\item[{\rm (ii)}]
$X=\{(x,y)|-1\leq x\leq 1,y=0\}\cup
\{(x,y)|x=0,-1\leq y\leq 1\}$.
\end{itemize}

\noindent
For a regular neighborhood $N(c_i)$ of $c_i$
in ${\bf R}^3$,
the pair $(N(c_i),D\cap N(c_i))$
is regarded as the image of
an immersion, say $\varphi$, of one of
the following manifold pairs:

\begin{itemize}
\setlength{\itemsep}{-3pt}
\item[{\rm (i)}]
$(B^2,X)\times [0,1]/
(x,0)\sim (x,1)$ for $x\in B^2$, and

\item[{\rm (ii)}]
$(B^2,X)\times [0,1]/
(x,0)\sim (-x,1)$ for $x\in B^2$.
\end{itemize}

\noindent
Let $c_i'$ be a loop or a pair of loops
immersed in $N(c_i)$
such that
$c_i'=\varphi\Bigl(\{z,-z\}\times [0,1]/\sim\Bigr)$
for some $z\in B^2\setminus X$.
We give orientations to $c_i$
and $c_i'$ such that $[c_i'] = 2[c_i] \in
H_1(N(c_i); {\bf Z})$.
See Figure $3$.
We put $C'=c_1'\cup\dots\cup c_m'$.
Since $C$ and $C'$ are mutually disjoint
1-cycles in ${\bf R}^3$,
the linking number Lk$(C,C')$
between $C$ and $C'$ is defined;
let $d$ be a 2-cycle in ${\bf R}^3$
with $\partial d =C$ such that $d$ and $C'$
intersect transversely.  Then Lk$(C, C')$
is the algebraic intersection  number of
$C'$ with $d$.
This number is well-defined modulo $4$;
it does not depend on a choice of
$z \in B^2 \setminus X$ and an orientation
of $c_i$ for each $i$.
Furthermore,
its congruence class modulo 4 remains unchanged under the 
Roseman moves. 
Hence, the mod $4$ reduction of Lk$(C,C')$
is an ambient isotopy invariant of
$F=F_1\cup F_2$.

\begin{figure}[htb]  \begin{center}
\mbox{
\epsfxsize=5in
\epsfbox{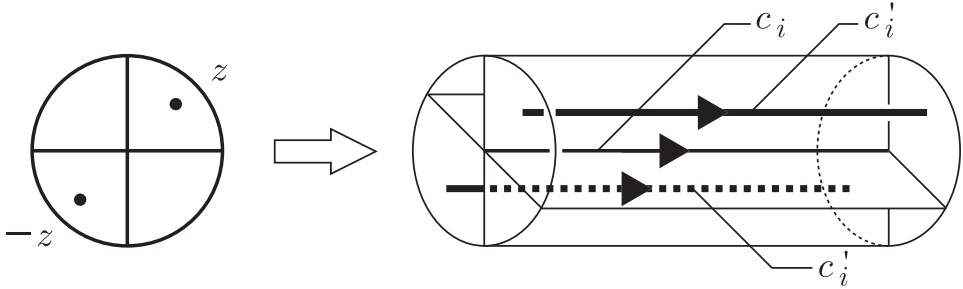} 
}
\end{center}
\caption{The double curve push-offs}
\end{figure}

\begin{definition}
{\rm The {\it double linking number}
between $F_1$ and $F_2$,
denoted by $d(F_1,F_2)$,
is a value in ${\bf Z}_4={\bf Z}/4{\bf Z}=\{0,1,2,3\}$
that is the linking number
${\rm Lk}(C,C')$ modulo $4$. }
\end{definition}

It is proved later that the double linking number
is asymmetric; $d(F_1,F_2)=-d(F_2,F_1)$.

Let $F=F_1\cup F_2\cup F_3$
be a $3$-component surface-link.
At a triple point on a surface diagram of $F$,
there are three sheets
called {\it top}, {\it middle}, and {\it bottom}
with respect to the projection direction.
A triple point is
{\it of type $(i,j,k)$}
if the top sheet comes from $F_i$,
the middle comes from $F_j$,
and the bottom comes from $F_k$,
where $i,j,k\in\{1,2,3\}$.
Let $N(i,j,k)$ denote
the number of the triple points of type $(i,j,k)$.
The mod $2$ reduction of $N(i,j,k)$
is preserved under Roseman moves and hence
an ambient isotopy invariant of $F$
provided $i\ne j$ and $j\ne k$
(possibly $i=k$), \cite{Satoh}.

\begin{definition}
{\rm The {\it triple linking number}
among $F_i,F_j$, and $F_k$,
denoted by $t(F_i,F_j,F_k)$,
is a value in ${\bf Z}_2={\bf Z}/2{\bf Z}=\{0,1\}$
that is the number $N(i,j,k)$ modulo $2$
provided $i\ne j$ and $j\ne k$.}
\end{definition}

\begin{lemma}
The ambient isotopy invariants
$e,d$ and $t$ are
unoriented bordism invariants.
\end{lemma}

{\it Proof.}
We take surface diagrams $D$ and $D'$
of surface-links $F$ and $F'$ respectively.
If $F$ and $F'$ are unorientedly bordant,
then  $D'$ is obtained from $D$
by a finite sequence of
moves from the following list:
\begin{itemize}
\setlength{\itemsep}{-3pt}
\item an ambient isotopy of ${\bf R}^3$,
\item a Roseman move,
\item adding or deleting embedded $2$-spheres
in ${\bf R}^3$ which are disjoint from $D$,
\item a $1$-handle surgery on $D$ in ${\bf R}^3$
whose core is a simple arc $\gamma$ with
$\gamma\cap D=\partial\gamma$
(see Figure $4$), and
\item a $2$-handle surgery on $D$ in ${\bf R}^3$
whose core is a simple $2$-disk $\delta$
with $\delta\cap D=\partial\delta$.
\end{itemize}

\noindent
Recall that Roseman moves do not change
the invariants $e,d$ and $t$.
The other deformations listed above
do not change the singularity set
of the diagram.
Hence $e,d$ and $t$ are unoriented bordism invariants.
\ \qed

\begin{figure}[htb]  \begin{center}
\mbox{
\epsfxsize=4.5in
\epsfbox{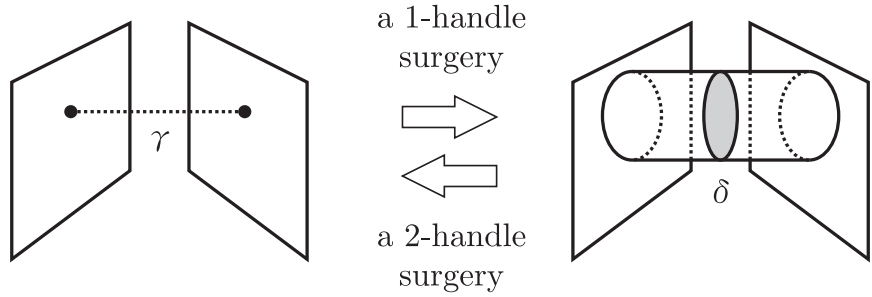} 
}
\end{center}
\caption{Attaching $1$ and $2$-handles}
\end{figure}


\section{1-Component Surface-Links}

A projective plane embedded
in ${\bf R}^4$ is {\it standard}
if it has a surface diagram
as shown in Figure $5$.
A non-orientable surface-knot is
said to be {\it trivial}
if it is a connected sum of some
standard projective planes in ${\bf R}^4$.
Two
trivial non-orientable surface-knots  $F$ and $F'$
are ambient isotopic if and only if
$e(F)=e(F')$ and $\chi(F)=\chi(F')$.
The following lemma is folklore.

\begin{figure}[htb]  \begin{center}
\mbox{
\epsfxsize=3in
\epsfbox{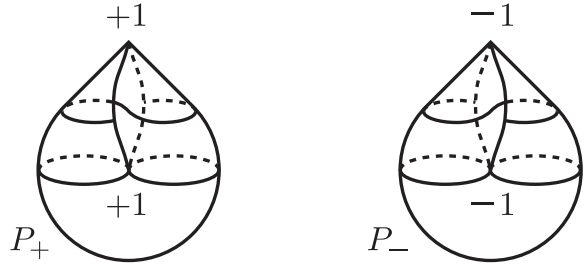} 
}
\end{center}
\caption{The positive and negative projective planes}
\end{figure}

\begin{lemma}
Two $1$-component surface-links
$F$ and $F'$ are unorientedly bordant
if and only if $e(F)=e(F')$.
\end{lemma}

{\it Proof.}
The only if part is obvious.
We prove the if part.
It is known that
any non-orientable surface-link
is transformed into
a trivial non-orientable surface-knot
by some $1$-handle surgeries,
\cite{Kamada}.
Thus we may assume that
$F$ and $F'$ are trivial non-orientable surface-knots
with $e(F)=e(F')$.
By Whitney's congruence,
we have $\chi(F)\equiv \chi(F')$ (mod~$2$).
Doing $1$-handle surgeries if necessary,
we may assume that $\chi(F)=\chi(F')$.
Then $F$ and $F'$ are ambient isotopic.
Thus original $1$-component surface-links
$F$ and $F'$ are unorientedly bordant.
\ \qed


\section{Necklaces}

We introduce a family of surface-links,
called {\it necklaces},
which is used to prove Theorem $1.2$.
In the upper $3$-space
${\bf R}^3_+=\{(x,y,z)|z\geq 0\}$,
we take a $3$-ball
$$B^3=\{(x,y,z)|x^2+y^2+(z-2)^2\leq 1\}.$$
Let $f=\{f_t\}_{0\leq t\leq 1}$
and $g=\{g_t\}_{0\leq t\leq 1}$
be ambient isotopies of $B^3$
which present
a $180^\circ$-rotation around the $z$-axis and
a $360^\circ$-rotation around the axis
$(y=0,z=2)$ respectively.
We put a Hopf link $k_1\cup k_2$ in $B^3$
as in Figure $6$
so that $f_1(k_i)=k_i$ 
and $g_1(k_i)=k_i$ for $i=1,2$.
By a {\it motion} of $k_1\cup k_2$,
we mean an ambient isotopy $h=\{h_t\}_{0\leq t\leq 1}$
of $B^3$ with $h_1(k_1\cup k_2)=k_1\cup k_2$.
Two motions $h$ and $h'$ are {\it equivalent}
if there is a $1$-parameter family
of motions between $h$ and $h'$.

\begin{figure}[htb]  \begin{center}
\mbox{
\epsfxsize=2in
\epsfbox{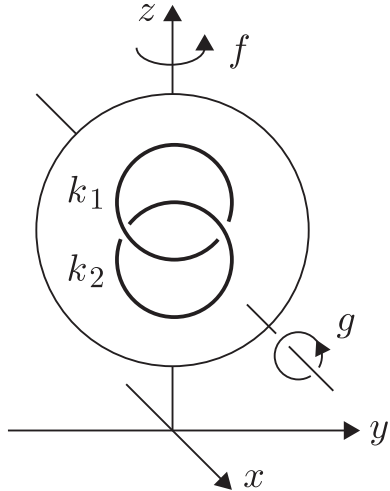} 
}
\end{center}
\caption{Spinning the Hopf link in two directions}
%
\end{figure}

We consider that ${\bf R}^4$ is obtained by
spinning ${\bf R}^3_+$ around $\partial{\bf R}^3_+$
by use of a map
$\mu:{\bf R}^3_+\times [0,1]\rightarrow {\bf R}^4$
defined by
$\Bigl((x,y,z) \times \{ t \}\Bigr)\mapsto (x,y,z\cos 2\pi t,z\sin 2\pi t)$.
  For integers $p,q$,
 we construct a $2$-component surface-link
$S^{p,q}=T_1\cup T_2$,
called a {\it strand}, as follows:
$$T_i=\mu\left(\bigcup_{t\in[0,1]}
\Bigl(f^p\cdot g^q\Bigr)_t(k_i)\times\{t\}\right)
\subset{\bf R}^4 \ (i=1,2).$$

\noindent
Each $T_i$ $(i=1,2)$ is homeomorphic to
a torus (or a Klein bottle, respectively)
if $p$ is even (or odd).

\begin{lemma}
{\rm (i)} Strands $S^{p,q}$ and $S^{p+2q,0}$ are
ambient isotopic.

{\rm (ii)} If $p\equiv p'$ {\rm (mod $4$)},
then two strands $S^{p,q}$ and $S^{p',q}$ are
ambient isotopic.
\end{lemma}

{\it Proof.}
By the belt trick (cf. \cite{Kauffman}),
the motion $g$ is equivalent to $f^2$, 
 and
$f^4$ is equivalent to the identity.
Hence we have the result.
\qed

\bigskip

The above lemma implies that ambient isotopy classes of strands are 
represented by $S^{p,q}$ with  $q=0$ and $p\in{\bf Z}_4$.
We shall abbreviate $S^{p,0}$ to $S^p$.

\begin{lemma}
For a strand $S^p=T_1\cup T_2$,
we have

{\rm (i)} $e(T_1)=e(T_2)=0\in{\bf Z}$,

{\rm (ii)} $d(T_1,T_2)=-d(T_2,T_1)=p
\in{\bf Z}_4$, and

{\rm (iii)} $t(T_1,T_2,T_1)=t(T_2,T_1,T_2)=p
\in{\bf Z}_2$.
\end{lemma}

{\it Proof.}
(i)
For each $i=1,2$,
we take a $2$-disk $D_i$ embedded in $B^3$ with
$\partial D_i=k_i$ and $f_1(D_i)=D_i$.
The image $\mu(D_i\times[0,1])$ is a $3$-manifold
whose boundary is $T_i$.
Thus $e(T_i)=0$.

(ii)
In Figure $7$,
we illustrate the motion $f$ of the Hopf link $k_1\cup k_2$.
Since we can obtain a diagram of $S^p$
by taking $p$ copies of the motion and connecting them,
we have $d(T_1,T_2)=p$.
Similarly, we have $d(T_2,T_1)=-p$.

(iii)
The motion in Figure $7$ contains
two Reidemeister moves of type III.
One of them corresponds to a triple point
of type (top, middle, bottom)$=(T_1,T_2,T_1)$ and
the other corresponds to that of $(T_2,T_1,T_2)$.
Thus we have $t(T_1,T_2,T_1)=t(T_2,T_1,T_2)=p$.
\qed

\bigskip

\begin{figure}[htb]  \begin{center}
\mbox{
\epsfxsize=5in
\epsfbox{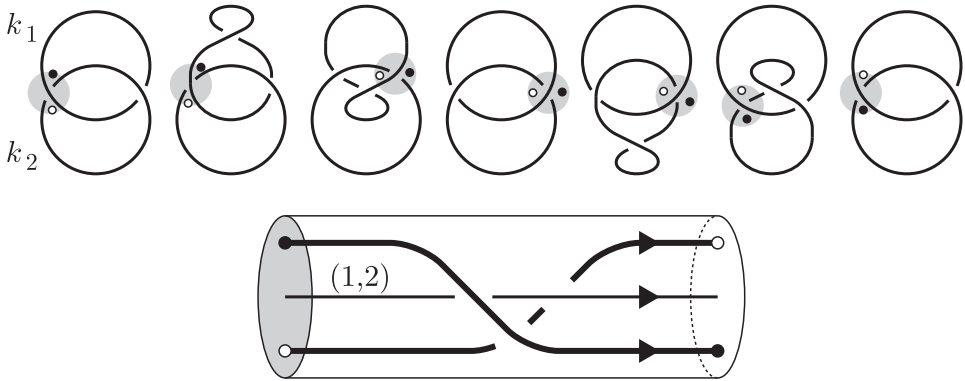} 
}
\end{center}
\caption{The movie of the strand $S^1$}
\end{figure}

For $m$ numbers $\{t_i\}_{i=1,\dots,m}$
with $0<t_1<\dots<t_m<1$,
we consider a surface-link 
defined 
as follows:
$$S^p\cup
\mu\Bigl(\partial B^3\times\{t_1\}\Bigr)
\cup\dots\cup
\mu\Bigl(\partial B^3\times\{t_m\}\Bigr).$$
We call such a surface-link a {\it necklace}
and a spherical component
$B_i=\mu\Bigl(\partial B^3\times\{t_i\}\Bigr)$
a {\it bead} of the necklace.

\begin{lemma}
Let $S^p\cup B_1\cup\dots\cup B_m$
be a necklace with the strand $S^p=T_1\cup T_2$.
For each $i=1,\dots,m$, we have
$$\left\{
\begin{array}{l}
t(T_2,T_1,B_i)=t(B_i,T_1,T_2)=1, \\
t(T_1,T_2,B_i)=t(B_i,T_2,T_1)=1, \\
t(T_1,B_i,T_2)=t(T_2,B_i,T_1)=0.
\end{array}\right.$$

\end{lemma}

{\it Proof.}
In a surface diagram,
a bead introduces four triple points
$\tau_1,\dots,\tau_4$
as shown in Figure $8$.
The top, middle and bottom sheets around $\tau_1$
come from $T_2$, $T_1$ and $B_i$ respectively.
For $\tau_2$, they are
$T_1$, $T_2$ and $B_i$.
For $\tau_3$, they are
$B_i$, $T_1$ and $T_2$.
For $\tau_4$, they are
$B_i$, $T_2$ and $T_1$.
\qed

\begin{figure}[htb]  \begin{center}
\mbox{
\epsfxsize=3in
\epsfbox{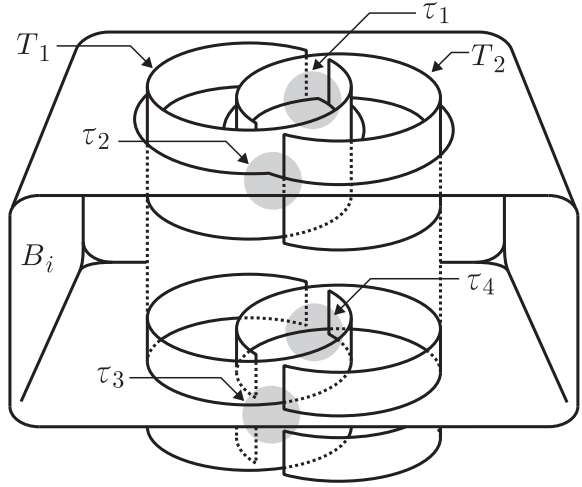} 
}
\end{center}
\caption{The local picture of a bead on a strand}
\end{figure}

\bigskip

We denote by
$N^p(i,j;k_1,\dots,k_m)$
an $n$-component surface-link which is a necklace
$S^p\cup B_1\cup\dots\cup B_m$
with the strand $S^p=T_1\cup T_2$ and
with a label
$\alpha:\{T_1,T_2,B_1,\dots,B_m\}\rightarrow
\{1,\dots,n\}$
such that
$\alpha(T_1)=i$, $\alpha(T_2)=j$, and
$\alpha(B_1)=k_1,\dots,\alpha(B_m)=k_m$.

\begin{lemma}
{\rm (i)}
$N^p(i,i;k_1,\dots,k_m)$
is unorientedly null-bordant.

{\rm (ii)}
$N^p(i,j;k_1,\dots,k_m)
\cong_A N^{-p}(j,i;k_1,\dots,k_m)$.

{\rm (iii)}
$N^p(i,j;k_1,k_2,k_3,\dots,k_m)
\simeq_B N^p(i,j;k_3,\dots,k_m)$
provided $k_1=k_2$.

{\rm (iv)}
$N^p(i,j;k_1,k_2,\dots,k_m)
\simeq_B N^{p+2}(i,j;k_2,\dots,k_m)$
provided $k_1=i$.

{\rm (v)}
$N^p(i,j;k_1,\dots,k_m)\amalg
N^{p'}(i,j;k_1',\dots,k_l')\simeq_B
N^{p+p'}(i,j;k_1,\dots,k_m,k_1',\dots,k_l')$.

{\rm (vi)}
$N^0(i,j;\emptyset)$ is
unorientedly null-bordant.

\end{lemma}

{\it Proof.}
(i)
Consider an annulus $A$ in int$B^3$
with $\partial A=k_1\cup k_2$
and $f_1(A)=A$.
The image $\mu(A\times[0,1])$ is
a $3$-manifold
whose boundary is the strand $S^p=T_1\cup T_2$.
Thus, we can remove $S^p$
and then all the beads $B_i$ $(i=1,\dots,m)$
up to unoriented bordism.

(ii)
By a $180^{\circ}$-rotation of $B^3$
around the axis $(y=0,z=2)$,
the components $k_1$ and $k_2$
are switched.
Then the direction of the rotation $f$
is reversed.

(iii)
We can remove the beads
$B_i=\mu\Bigl(\partial B^3\times\{t_i\}\Bigr)$
$(i=1,2)$ from $F$ by a $3$-manifold
$\mu\Bigl(\partial B^3\times[t_1,t_2]\Bigr)$.

(iv)
We perform a $1$-handle surgery
between $K_1$ and $B_1$
as shown in Figure $9(1)\rightarrow(2)$.
The diagram illustrated in $(2)$ is
ambient isotopic to that in (3).
This surface is realized by a motion as in $(4)$,
which is the $360^\circ$-rotation around the axis
$(x=0,z=2)$.
This motion is equivalent to $f^2$.

(v)
Connect the strands of
$N^p(i,j;k_1,\dots,k_m)$ and
$N^{p'}(i,j;k_1',\dots,k_l')$
up to unoriented bordism
to obtain
$N^{p+p'}(i,j;k_1,\dots,k_m,k_1',\dots,k_l')$
(see Lemma $2.3$ of \cite{CKSS2}).

(vi)
Consider the Hopf link $k_1\cup k_2$
to use the definition of a strand.
Let $\gamma_i$ $(i=1,2)$
be disjoint simple arcs
in ${\bf R}_3^+$ connecting between
$k_i$ and $\partial{\bf R}^3_+$
in an obvious way.
Then the $2$-handle surgeries along the $2$-disks
$\mu\Bigl(\gamma_1\times[0,1]\Bigr)$ and
$\mu\Bigl(\gamma_2\times[0,1]\Bigr)$ make $S^0$
a split union of $2$-spheres.
\qed

\begin{figure}[htb]  \begin{center}
\mbox{
\epsfxsize=5in
\epsfbox{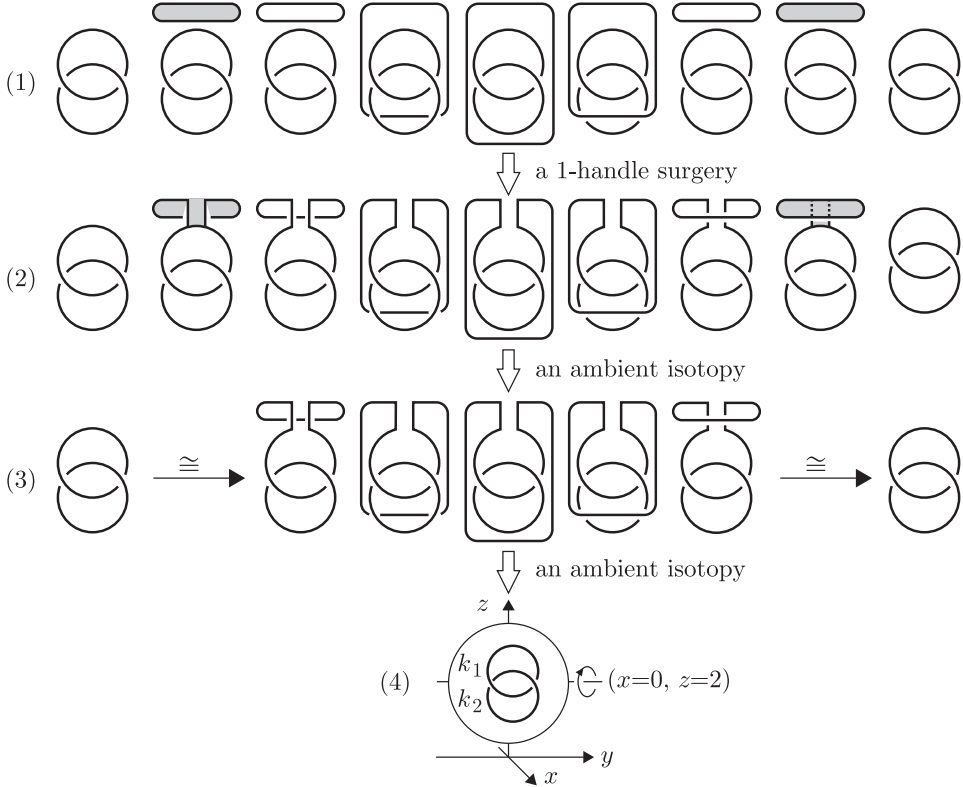} 
}
\end{center}
\caption{The movies of surgeries and isotopies}
\end{figure}


\section{Crossing Changes}

In this section,
we study a crossing change
along a double curve of a surface diagram.
The idea is similar to a crossing change for a classical link
in ${\bf R}^3$.
For a classical link,
a crossing change can be
realized by $1$-handle(=band) surgeries
as shown in Figure $10$(i).
Hence, any classical link $k_1\cup\dots\cup k_n$
is bordant to a split union
$k_1\amalg\dots\amalg k_n\amalg($Hopf links$)$.

\begin{figure}[htb]  \begin{center}
\mbox{
\epsfxsize=5in
\epsfbox{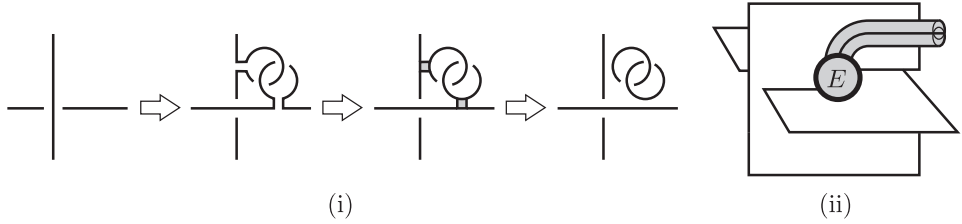} 
}
\end{center}
\caption{A movie of a crossing change and the resulting surface}
\end{figure}

\begin{lemma}
Any $n$-component surface-link
$F_1\cup\dots\cup F_n$ is
unorientedly bordant to a split union
$F_1\amalg\dots\amalg F_n\amalg({\rm necklaces})$.
\end{lemma}

In a surface diagram, we use a symbol $E$
(as in Figure $10$(ii)) for a
local diagram which is
obtained from Figure $10$(i) by
regarding it a motion picture.

Let $c$ be a double curve of a surface diagram
such that it is an immersed loop in ${\bf R}^3$.
As mentioned in Section $2$,
a regular neighborhood $N(c)$ of $c$
is regarded as the image of an immersion $\varphi$ of
  $B^2 \times [0,1]/ \sim$.
See Figure $11$.
For each $c$,
we choose a pair of diagonal regions
$Y$ of $B^2\setminus X$
and put
$R(c)=\varphi(Y\times[0,1]/\sim)$.

\begin{figure}[htb]  \begin{center}
\mbox{
\epsfxsize=3in
\epsfbox{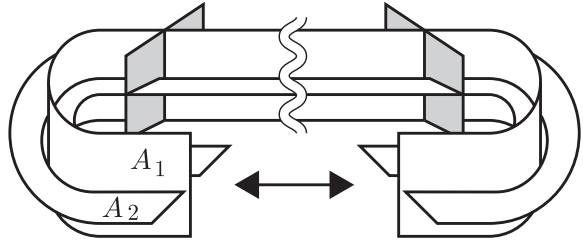} 
}
\end{center}
\caption{A neighborhood of a  double curve}
\end{figure}

Let $A_1$ and $A_2$ be the sheets
which intersects along the curve $c$
such that $A_1$ is higher than $A_2$
with respect to the projection direction;
that is, $c$ is of type (upper, lower)=$(A_1,A_2)$.
To prove Lemma $5.1$,
we introduce four kinds of local deformations
under which unoriented bordism classes are preserved.

{\bf (1) Local crossing change.}
Consider two $1$-handle surgeries
as shown in Figure $12$
in the motion picture method.
We call this a {\it local crossing change} along $c$.
We always assume that
the `local strand' (Hopf link)$\times[0,1]$
obtained by a local crossing change
is in the specified region $R(c)$.
Let $t$ be a triple point on $c$ and
$H$ a sheet which is transverse to $c$ at $t$.
If $H$ is a top or bottom sheet,
then we can perform a local crossing change
along the curve $c$ as in Figure $13$.
For example,
if $H$ is a top sheet,
then the type of $t$ is changed from
(top, middle, bottom)=$(H,A_1,A_2)$ to
$(H,A_2,A_1)$, and the local strand
goes under the sheet $H$.

\begin{figure}[htb]  \begin{center}
\mbox{
\epsfxsize=5in
\epsfbox{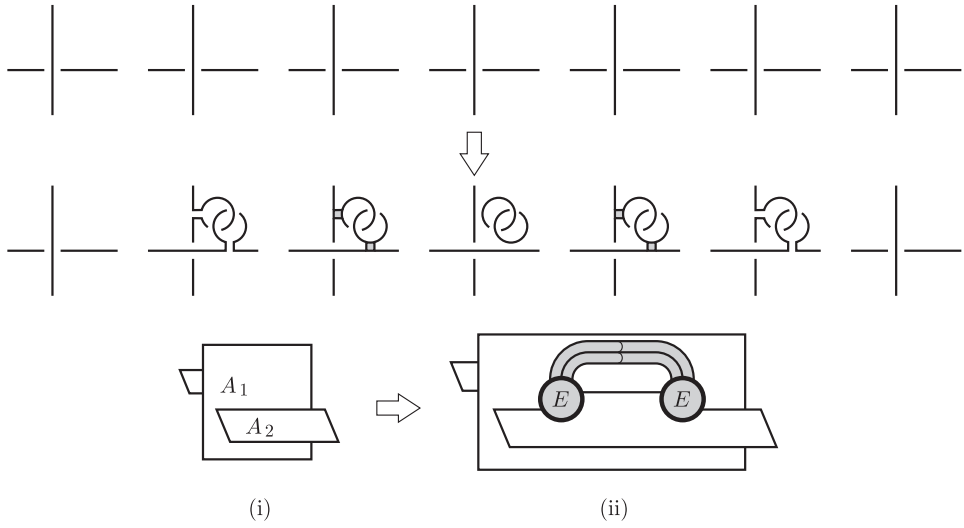} 
}
\end{center}
\caption{The local crossing change}
\end{figure}

\begin{figure}[htb]  \begin{center}
\mbox{
\epsfxsize=4in
\epsfbox{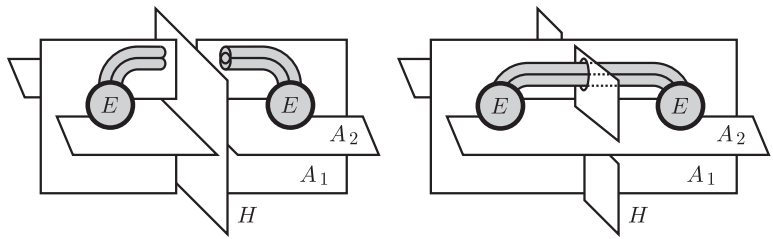} 
}
\end{center}
\caption{A local crossing change near a triple point}
\end{figure}

\clearpage

{\bf (2) End-change.}
Consider a composite of three ambient isotopies
as shown in Figure $14$.
We call this an {\it end-change} of a local strand.
This deformation moves
an end-part of the local strand
into the diagonal region of $R(c)$.
The new local strand has
additional intersections with $A_1$ and $A_2$.
In the bottom of Figure $14$,
the boxed ``$f$'' means a local diagram
corresponding to a $180^\circ$-rotation
of the Hopf link.

\begin{figure}[htb]  \begin{center}
\mbox{
\epsfxsize=5in
\epsfbox{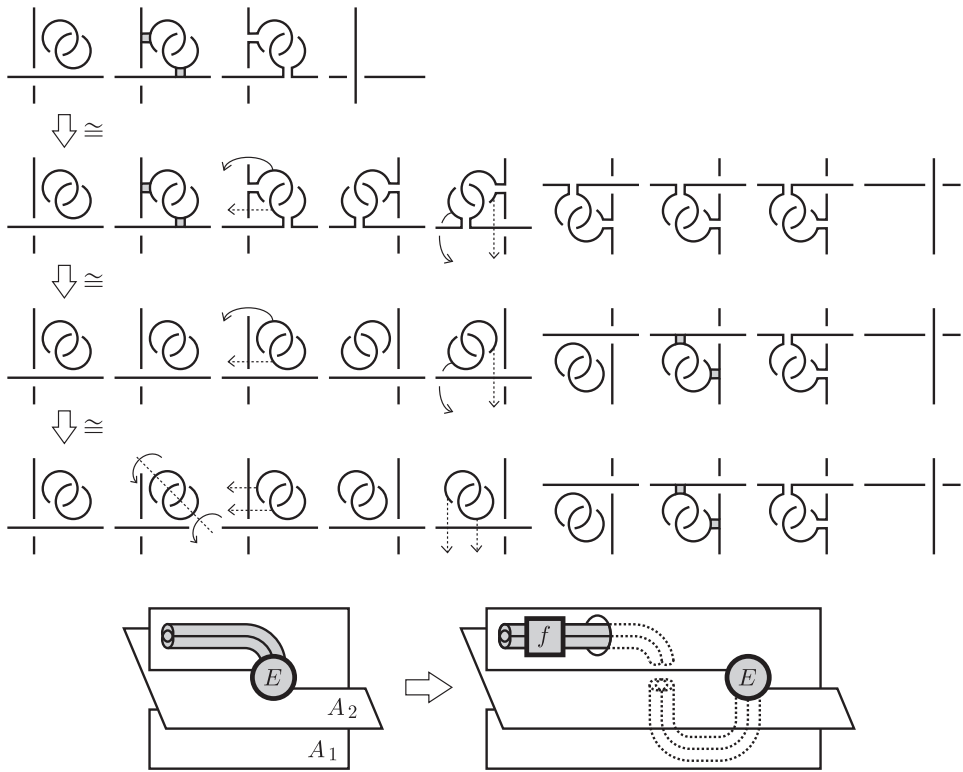} 
}
\end{center}
\caption{The movies and diagrams determining an end change}
\end{figure}

{\bf (3) Canceling adjacent ends.}
Assume that two adjacent ends
of local strands are in the same region of $R(c)$.
Consider a deformation illustrated in Figure $15$,
which is realized by two $2$-handle surgeries
on $A_1$ and $A_2$.
We call this a {\it canceling of adjacent ends}
of local strands along $c$.

\begin{figure}[htb]  \begin{center}
\mbox{
\epsfxsize=4in
\epsfbox{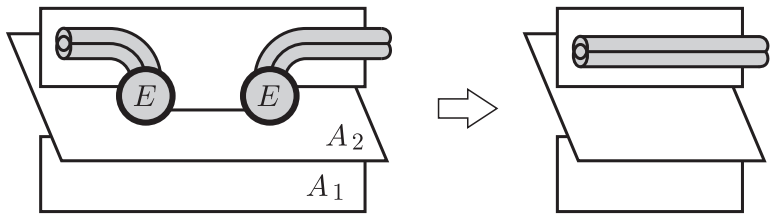} 
}
\end{center}
\caption{Canceling adjacent ends for a crossing change}
\end{figure}

{\bf (4) Making a bead.}
Assume that a strand intersects a sheet $H$
transversely.
In case that the strand is under $H$,
we consider a $2$-handle surgery
on $H$ as shown in Figure $16$.
This deformation makes the strand
over $H$ with producing a bead.

\begin{figure}[htb]  \begin{center}
\mbox{
\epsfxsize=4in
\epsfbox{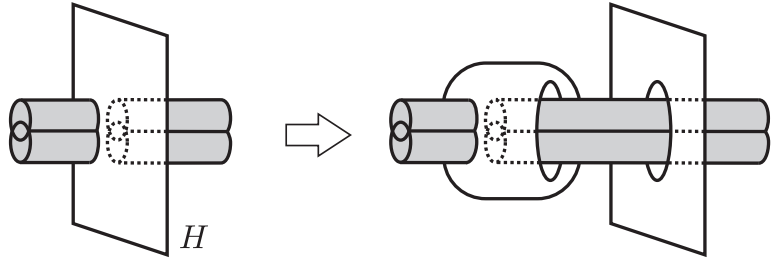} 
}
\end{center}
\caption{Making a bead}
\end{figure}

In a surface diagram of
an $n$-component surface-link
$F=F_1\cup\dots\cup F_n$,
a double curve of type $(i,j)$
is {\it preferred} if $i\leq j$,
and a triple point of type $(i,j,k)$ is {\it preferred}
if $i\leq j$ and $j\leq k$.

\bigskip

{\it Proof of Lemma $5.1$.}
Consider a surface diagram of $F=F_1\cup\dots\cup F_n$.
If there is a non-preferred double curve
without triple points
(that is an embedded loop),
we apply a local crossing change along the curve
followed by canceling adjacent ends of the local strand.
This makes the double curve preferred and
yields a strand (necklace without beads),
which is separated from $F$.
Thus we may assume that
there is at least one triple point on
any non-preferred double curve of the surface diagram.

For each non-preferred triple point $t$
of type $(i,j,k)$,
we perform local crossing changes so that
$t$ changes into a preferred triple point
as follows:

\begin{itemize}
\setlength{\itemsep}{-3pt}
\item
If $j<i\leq k$, then
we perform a local crossing change
along the double curve of type $(i,j)$
over the bottom sheet labeled $k$.
\item
If $i\leq k<j$, then
we perform a local crossing change
along the double curve of type $(j,k)$
under the top sheet labeled $i$.
\item
If $k<i\leq j$, then
we perform a local crossing change
along the double curve of type $(j,k)$
and then perform another
along the curve of type $(i,k)$.
\item
If $j\leq k<i$, then
we perform a local crossing change
along the double curve of type $(i,j)$
and then perform anther
along the curve of type $(i,k)$.
\item
If $k<j<i$, then
we perform a local crossing change
along the double curve of type $(i,j)$,
next along the curve of type $(i,k)$,
and then along the curve of type $(j,k)$.
\end{itemize}

\noindent
Figure $17$ shows the case of $k<j<i$.

\begin{figure}[htb]  \begin{center}
\mbox{
\epsfxsize=5in
\epsfbox{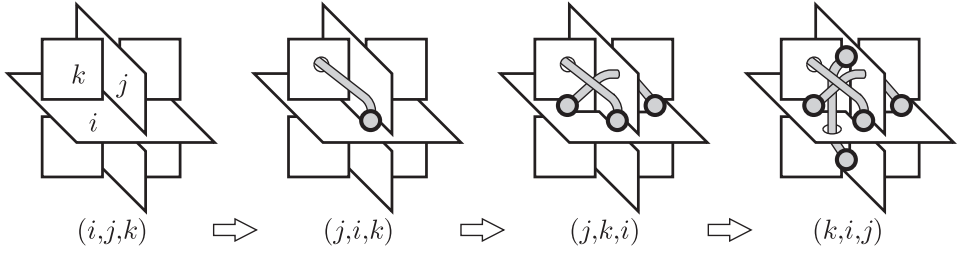} 
}
\end{center}
\caption{A three-fold crossing change at a triple point of type $k<j< i$}
\end{figure}

All the local crossing changes 
that are 
described above
are 
performed 
along non-preferred double curves.
After applying suitable end-changes,
we can cancel all the adjacent ends of local strands
along non-preferred double curves.
Then we obtain a surface diagram of
$(F_1\amalg\dots\amalg F_n)\cup($strands$)$;
for any double curves between $F_i$ and $F_j$ are preferred.
By making beads if necessary,
we can split necklaces from $F_1\amalg\dots\amalg F_n$.
Thus $F$ is unorientedly bordant to
$F_1\amalg\dots\amalg F_n\amalg($necklaces$)$.
\qed

\bigskip

\begin{lemma}
$N^0(i,j;k)\simeq_B N^0(k,i;j)\amalg N^0(k,j;i)$.
\end{lemma}

{\it Proof.}
Consider a surface diagram of $N^0(i,j;k)$
which contains the local diagram
illustrated in Figure $8$.
In the local diagram,
Along the double curves of type $(i,k)$ and $(j,k)$,
we perform (global) crossing changes
as in the proof of Lemma $5.1$.
Then $N^0(i,j;k)$ is unorientedly bordant to
a split union of $N^0(i,k;j)\cong_A N^0(k,i;j)$
along the curve of type $(i,k)$,
$N^0(j,k;i)\cong_A N^0(k,j;i)$ along the curve of
type $(j,k)$, and a surface-link $F$
obtained by the crossing change.
Then $F$ is a split union of $S^0(i,j)$
and a trivial $2$-sphere labeled $k$,
which is unorientedly null-bordant
by Lemma $4.4$(vi).
\qed


\section{Proof of Theorem 1.2}

Let $P_+$ and $P_-$
be the standard projective plane in ${\bf R}^4$
with $e(P_+)=2$ and $e(P_-)=-2$ respectively.
See Figure $5$.
We denote by $P^m$
the connected sum of $m$ copies of
$P_+$ if $m>0$,
$-m$ copies of $P_-$ if $m<0$,
and the empty if $m=0$.
Regarding $P^m$ as an $n$-component surface-link
with a label $\alpha(P^m)=i$,
we denote it by $P^m(i)$.
Regarding a strand $S^p=T_1\cup T_2$
as an $n$-component surface-link
such that $\alpha(T_1)=i$ and $\alpha(T_2)=j$,
we denote it by $S^p(i,j)$.
Then $S^p(i,j)=N^p(i,j;\emptyset)$
in the notation used in Lemma $4.4$.

\begin{lemma}
Any $n$-component surface-link $F$ is
unorientedly bordant to a split union
$$\Biggl(\coprod_{i\in\Gamma_1} P^{m_i}(i)\Biggr)\coprod
\Biggl(\coprod_{(i,j)\in\Gamma_2} S^{p_{ij}}(i,j)\Biggr)\coprod
\Biggl(\coprod_{(i,j,k)\in\Gamma_3} N^0(i,j;k)\Biggr).$$
Here,
$\Gamma_i$ $(i=1,2,3)$
is a subset of the $i$-fold Cartesian product
of $\{1,\dots,n\}$,
$m_i\in{\bf Z}$ $(i\in\Gamma_1)$ and
$p_{ij}\in{\bf Z}_4$ $\Bigl((i,j)\in\Gamma_2\Bigr)$
satisfying the following.

\begin{itemize}
\setlength{\itemsep}{-3pt}
\item[{\rm (i)}]
$m_i\ne 0$ for any $i\in\Gamma_1$.

\item[{\rm (ii)}]
$i<j$ and $p_{ij}\ne 0$
for any $(i,j)\in\Gamma_2$.

\item[{\rm (iii)}]
$i<j<k$ or $i<k<j$ for any $(i,j,k)\in\Gamma_3$.
\end{itemize}
\end{lemma}

{\it Proof.}
By Lemma $5.1$,
any $n$-component surface-link $F=F_1\cup\dots\cup F_n$
is unorientedly bordant to
$(F_1\amalg\dots\amalg F_n)\amalg F'$,
where $F'$ is a split union of necklaces.
Put
$\Gamma_1=\{i\ |e(F_i)\ne 0\}$ and  $m_i=e(F_i)/2\in{\bf Z}$
$(i\in\Gamma_1)$.
By Lemma $3.1$,
we see that
$F_1\amalg\dots\amalg F_n$ is unorientedly bordant to
$\coprod_{i\in\Gamma_1} P^{m_i}(i)$
satisfying the condition (i).

Since a necklace
$N^p(i,j;k_1,\dots,k_m)$ is unorientedly bordant to
$S^p(i,j)\amalg N^0(i,j;k_1)\amalg\dots
\amalg N^0(i,j;k_m)$
by Lemma $4.4$(v),
$F'$ is unorientedly bordant to $F''\amalg F'''$
such that $F''$ is a split union of
some $S^p(i,j)$'s and
$F'''$ is a split union of some $N^0(i,j;k)$'s.

We may assume that $i<j$
for any $S^p(i,j)$ appearing in $F''$
by Lemma $4.4$(i) and (ii).
Moreover,
by Lemma $4.4$(v) and (vi),
we see that there exist a subset
$\Gamma_2\subset\{1,\dots,n\}^2$ and
$p_{ij}\in{\bf Z}_4$
such that $F''$ is unorientedly bordant to
$\coprod_{(i,j)\in\Gamma_2} S^{p_{ij}}(i,j)$
and that the condition (ii) is satisfied.

By Lemma $4.4$(i), (ii), (iv) and (vi),
we may assume that $i<j$,
$i\ne k$ and $j\ne k$ for any $N^0(i,j;k)$
appearing in $F'''$.
Applying Lemma $5.2$ for $N^0(i,j;k)$
with $k<i<j$,
we may assume that $(i,j,k)$ satisfies
the condition (iii) for any $N^0(i,j;k)$
appearing in $F'''$.
By Lemma $4.4$(iii), (v) and (vi),
we see that there exists a subset
$\Gamma_3\subset\{1,\dots,n\}^3$
such that $F'''$ is unorientedly bordant to
$\coprod_{(i,j,k)\in\Gamma_3} N^0(i,j;k)$
and that the condition (iii) is satisfied.
\qed

\bigskip

For the unoriented bordism group $UL_{4,n}$,
we consider three types of
homomorphisms $e_i$ $(i=1,\dots,n)$,
$d_{ij}$ $(i\ne j)$, and
$t_{ijk}$ ($i\ne j$ and $j\ne k$) as follows:
$$\left\{\begin{array}{rll}
e_i:&UL_{4,n}\longrightarrow{\bf Z}
&{\rm for}\ [F]\mapsto e(F_i)/2,\\
d_{ij}:&UL_{4,n}\longrightarrow{\bf Z}_4
&{\rm for}\ [F]\mapsto d(F_i,F_j), \\
t_{ijk}:&UL_{4,n}\longrightarrow{\bf Z}_2
&{\rm for}\ [F]\mapsto t(F_i,F_j,F_k),
\end{array}\right.
$$
where $F=F_1\cup\dots\cup F_n$.

\begin{lemma}
For an $n$-component surface-link $F$,
let $\Gamma_i$ $(i=1,2,3)$,
$m_i\in{\bf Z}$ $(i\in\Gamma_1)$ and
$p_{ij}\in{\bf Z}_4$ $\Bigl((i,j)\in\Gamma_2\Bigr)$
be as in Lemma $6.1$.
Then we have the following.
\begin{itemize}
\setlength{\itemsep}{-2pt}
\item[{\rm (i)}]
$e_i([F])=m_i$ if $i\in\Gamma_1$ and
$e_i([F])=0$ if $i \not\in\Gamma_1$.
\item[{\rm (ii)}]
For $i<j$,
$d_{ij}([F])=p_{ij}$ if $(i,j)\in\Gamma_2$
and $d_{ij}([F])=0$ if $(i,j)\not\in\Gamma_2$.

\item[{\rm (iii)}]
For $i<j<k$,
$t_{ijk}([F])=1$ if $(i,j,k)\in\Gamma_3$
and $t_{ijk}([F])=0$ if $(i,j,k)\not\in\Gamma_3$.

\item[{\rm (iv)}]
For $i<k<j$,
$t_{ijk}([F])=1$ if $(i,j,k)\in\Gamma_3$
and $t_{ijk}([F])=0$ if $(i,j,k)\not\in\Gamma_3$.

\end{itemize}
\end{lemma}

{\it Proof.}
(i)
Since $e(P^m)=2m$,
we have $e_i([F])=m_i$ if $i\in\Gamma_1$
and otherwise $e_i([F])=0$.

(ii)
This follows from Lemma $4.2$(ii).

(iii) and (iv)
Note that $(j,i,k)$, $(j,k,i)$,
$(k,i,j)$, $(k,j,i)\notin\Gamma_3$.
Since $t_{ijk}\Bigl([N^0(i,j;k)]\Bigr)=1$
and $t_{ijk}\Bigl([N^0(i,k;j)]\Bigr)=0$
by Lemma $4.3$,
$t_{ijk}([F])=1$ if and only if
$(i,j,k)\in\Gamma_3$.
\qed

\bigskip

{\it Proof of Theorem $1.2$.}
Consider a homomorphism
$$UH:UL_{4,n}\longrightarrow
(\underbrace{{\bf Z}\oplus\dots
\oplus{\bf Z}}_{n})
\oplus
(\underbrace{{\bf Z}_4\oplus\dots
\oplus{\bf Z}_4}_{\frac{n(n-1)}{2}})
\oplus
(\underbrace{{\bf Z}_2 \oplus\dots
\oplus{\bf Z}_2}_{\frac{n(n-1)(n-2)}{3}})$$
defined by
$UH=\Bigl(\bigoplus_{i=1}^n e_i\Bigr)\oplus
\Bigl(\bigoplus_{i<j} d_{ij}\Bigr)\oplus
\Bigl(\bigoplus_{i<j<k\mbox{ or }i<k<j}
t_{ijk}\Bigr)$.
This homomorphism is injective by Lemma $6.2$.
Also,
$UH$ is surjective;
indeed,
$UH\Bigl([P^1(i)]\Bigr)$ $(i=1,\dots,n)$,
$UH\Bigl([S^1(i,j)]\Bigr)$ $(i<j)$, and
$UH\Bigl([N^0(i,j;k)]\Bigr)$
$(i<j<k$ or $i<k<j)$ are
generators of
${\bf Z}$, ${\bf Z}_4$, and ${\bf Z}_2$
respectively.
\qed

\bigskip

In the definition of the homomorphism $UH$,
we do not use all double linking numbers and
triple linking numbers.
The unused ones are determined as follows.

\begin{proposition}
For distinct $i,j,k$ and
an $n$-component surface-link $F$,
we have

{\rm (i)}
$d_{ji}([F])=-d_{ij}([F])$,

{\rm (ii)}
$t_{iji}([F])=t_{jij}([F])=
\lambda\Bigl(d_{ij}([F])\Bigr)$,
where $\lambda:{\bf Z}_4\rightarrow{\bf Z}_2$
is the natural projection,

{\rm (iii)}
$t_{ijk}([F])=t_{kji}([F])$, and

{\rm (iv)}
$t_{jik}([F])+t_{ijk}([F])+t_{ikj}([F])=0$.

\end{proposition}

{\it Proof.}
It is sufficient to prove (i) and (ii)
in case that $F$ is as in Lemma $6.1$.
We use Lemma $4.2$.
We have $t_{iji}([F])=p_{ij}\in{\bf Z}_2$ if $i<j$, and
$t_{iji}([F])=p_{ji}\in{\bf Z}_2$ if $i>j$.
On the other hand,
we have $d_{ij}([F])=p_{ij}\in{\bf Z}_4$ if $i<j$,
and
$d_{ij}([F])=-p_{ji}\in{\bf Z}_4$ if $i>j$.
Hence, we have
$\lambda\Bigl(d_{ij}([F])\Bigr)=t_{iji}([F])$.
Similarly,
since $d_{ij}([F])=p_{ij}\in{\bf Z}_4$ and
$d_{ji}([F])=-p_{ij}\in{\bf Z}_4$ for $i<j$,
we have $d_{ij}([F])=-d_{ji}([F])$.
(iii) and (iv) are proved
in Theorem $3.2$ of \cite{Satoh}.
\qed

\bigskip

We consider the homomorphism
$f:L_{4,n}\rightarrow UL_{4,n}$
induced by the map
forgetting the orientations of surface-links.
For an oriented $n$-component surface-link $F$,
we can define two kinds of bordism invariants;
double linking invariants
$D_{ij}:L_{4,n}\rightarrow{\bf Z}_2={\bf Z}/2{\bf Z}$
and triple linking invariants
$T_{ijk}:L_{4,n}\rightarrow{\bf Z}$
(cf. \cite{CKSS2}).
Then Sanderson's isomorphism
$$H:L_{4,n}\longrightarrow
(\underbrace{{\bf Z}_2\oplus\dots
\oplus{\bf Z}_2}_{\frac{n(n-1)}{2}})
\oplus
(\underbrace{{\bf Z} \oplus\dots
\oplus{\bf Z}}_{\frac{n(n-1)(n-2)}{3}})$$
is given by
$H=\Bigl(\bigoplus_{i<j} D_{ij}\Bigr)\oplus
\Bigl(\bigoplus_{i<j<k\mbox{ or }i<k<j}
T_{ijk}\Bigr)$.
 From the definitions of these invariants,
the forgetful map $f$ is regarded as
$$(\oplus 0)\oplus(\oplus\kappa)\oplus
(\oplus\nu):(\oplus\{0\})\oplus(\oplus{\bf Z}_2)
\oplus(\oplus{\bf Z})\rightarrow
(\oplus{\bf Z})\oplus(\oplus{\bf Z}_4)
\oplus(\oplus{\bf Z}_2),$$
under the isomorphisms $U$ and $UH$,
where $\kappa:{\bf Z}_2\rightarrow{\bf Z}_4$
is the natural inclusion and
$\nu:{\bf Z}\rightarrow{\bf Z}_2$ is
the natural projection.

\bigskip

\noindent
{\large\bf Acknowledgments} \
JSC is being supported by NSF grant DMS-9988107.
SK is being supported by Fellowships
from the Japan Society for the Promotion of Science.
MS is being supported by NSF grant DMS-9988101.
The authors would like to thank Uwe Kaiser for providing helpful information. 


\newcommand{\bysame}{%
\leavevmode\hbox to 3em{\hrulefill}\,}

\end{document}